\documentclass[leqno,12pt]{article}

\usepackage[utf8]{inputenc}
\usepackage[T1]{fontenc}
\usepackage[english]{babel}
\usepackage[intlimits]{amsmath}
\usepackage{amsfonts}
\usepackage{amssymb}
\usepackage{bbm}

\usepackage{textcomp}
\usepackage{array}
\usepackage{enumerate}

\usepackage[amsmath, thmmarks]{ntheorem}
\usepackage{multicol}
\usepackage{units}
\usepackage{nicefrac}

\newcounter{theoremcounter}

\theoremstyle{plain}
\theoremseparator{.}

\theoremstyle{plain}
\newtheorem{lemma}[theoremcounter]{Lemma}
\newtheorem{theorem}[theoremcounter]{Theorem}

\theorembodyfont{\normalfont}
\theoremsymbol{\ensuremath{\Box}}
\newtheorem*{eproof}{Proof}

\newcommand{\N}{\mathbb{N}}
\newcommand{\Z}{\mathbb{Z}}
\newcommand{\Q}{\mathbb{Q}}

\newcommand{\tr}{\mathrm{tr}}

\newcommand{\Mat}[1]{\mathrm{M}_{#1}}
\newcommand{\GL}[1]{\mathrm{GL}_{#1}}

\newcommand{\RGL}[1]{\widetilde{\mathrm{GL}}_{#1}}
\newcommand{\Sp}[1]{\mathrm{Sp}_{#1}}

\newcommand{\Inv}[1]{\mathrm{Inv}_{#1}}
\newcommand{\RInv}[1]{\mathrm{Inv}_{#1} ^{\mathrm{Res}}}

\newcommand{\compl}[2]{#1_{\mathfrak{#2}}}
\newcommand{\ccompl}[2]{#1_{\bar{\mathfrak{#2}}}}
\newcommand{\adelic}[1]{#1_{\mathbb{A}}}

\begin{document}
\title{Hecke algebras related to the unimodular and modular groups over hermitian fields and definite quaternion algebras}
\author{Martin Raum\vspace{0.5em} \\Lehrstuhl A f\"ur
Mathematik\\RWTH Aachen University, Germany}
\date{\today}
\maketitle

\begin{abstract}
We investigate the structure of the Hecke algebras related to the unimodular and modular group over hermitian fields and definite quaternion algebras. In particular we show that in general there is no decomposition into primary components.
We can give a set of generators and in some special cases we deduce the law of interchange of the Siegel $\Phi$-operator.
\end{abstract}
\section{Introduction}
Many arithmetical informations, i.e. the attached $\zeta$-functions, of modular forms are related to the action of the associated Hecke algebras. These Hecke algebras have been investigated for many principle ideal domains \cite{Kr90, Kr87_ha}. Nevertheless the case of non PIDs remains mysterious until now. In \cite{Eb08} evidence has been provided that it behaves differently.

We prove that the Hecke algebras in question for hermitian fields and quaternionic algebras are subalgebras of their adelic analogs. The point of view suggested by this discovering makes it possible to deduce a set of generators and determine how far the algebras are from having a primary decomposition.

In the general PID case, in addition, we are able to deduce the yet unproven law of interchange for the Siegel $\Phi$-operator.

In \cite{R09} we proved that the equivalence class of any unimodular or modular matrix in quadratic extensions of fields and quaternion algebras is completely determined by its adelic elementary divisors. In addition we investigated which adelic elementary divisors may occur. This statement can be used to investigate the structure of the related Hecke algebras.

\section{Embeddings of integral Hecke algebras}\label{sec:subheckealgebra}

Let $K$ be a number field and denote its ring of integers $\mathfrak{o}_K$.
Depending on the case we treat, $\Omega$ will either denote a field extension of $K$ or a central simple $K$-algebra and $\Lambda$ will always be a maximal order of $\Omega$. Note that we don't exclude $K = \Omega$.
The set of finite places of $\Omega$ will be denoted $\adelic{\Omega}^f$ and the set of all infinite places $\adelic{\Omega}^\infty$. For $\mathfrak{p} \in \adelic{\Omega}^f$ the completions of $\Omega,\,\Lambda$ etc. with respect to $\mathfrak{p}$ will be denoted $\compl{\Omega}{p},\,\compl{\Lambda}{p}$ etc.. The $\adelic{\Omega}^\infty$-integral adelic algebras associated to $\Omega$ and $\Lambda$ will be denoted $\adelic{\Omega}$ and $\adelic{\Lambda}$.

Let $n \in \N$ and whenever $\Lambda$ is non commutative assume $n \ge 2$.
Define $\RGL{n} (\adelic{\Lambda}) \subseteq \GL{n} (\adelic{\Lambda})$ to be the group generated by all matrices $I_n + a I^{(n,n)} _{i,j}$ where $a \in \adelic{\Lambda}$ and $i \ne j$.
Here $I_n$ is the identity matrix and $I^{(n,n)}_{i,j}$ is the $n\times n$ matrix with only zeros but a single one in the $i$-th row and $j$-th column.
Moreover let $\Inv{n} (\adelic{\Lambda})$ be the monoid of all matrices which are invertible over the fraction ring of $\adelic{\Lambda}$.
In analogy with this case define $\RGL{n} (\compl{\Lambda}{p}),\,\Inv{n} (\compl{\Lambda}{p})$ etc..

Two elements $a,\,b \in \Lambda$ are called similar, denoted $a \sim b$, if there are units $u,\,v \in \Lambda^\times$ such that $a = u b v$.

If there is an involution $\bar \cdot$ of the algebra extension $\Omega | K$ which fixes $K$ elementwise and preserves maximal orders we can define the set of all matrices with similitude $m \in K^\times$.
To do so, denote the matrix transpose by $M^\tr$ and define $M^* := \overline{M} ^\tr$ for any matrix $M \in \Mat{n} (\Omega)$.
\begin{align*}
  \Delta_n (\Omega, m)
:=
  \{ M \in \Mat{2 n} (\Omega) : M^* J M = m J\}
\quad
\text{with}
\quad
  J = \left(\begin{matrix}
      0^{(n,n)} & I_n \\
      - I_n     & 0^{(n,n)}
      \end{matrix}\right)
.
\end{align*}
The monoid $\Delta_n (\Omega) := \bigcup_{m \in \mathfrak{o}_K \setminus \{0\}} \Delta_n (\Omega, m)$ of all matrices with integral similitude contains the modular group $\Sp{n} (\Omega) := \bigcup_{m \in \mathfrak{o}_K ^\times} \Delta_n(\Omega, m)$.
We point out that most authors let $K = \Q$ and only treat similitudes $m \in \N$. This is no restriction, as we will see, since using $\Delta_n (\Z, 1)$ every matrix may be reduced to a diagonal matrix. Hence it is possible to choose a multiplicative set of representatives of $(\mathfrak{o}_K \setminus \{0\}) / \mathfrak{o}_K ^\times$ and for $K = \Q$ this yields the aforementioned choice.

Recall the definition of a Hecke algebra. Let $G$ be a group and let $H \subseteq G$ be a monoid such that each double coset $G h G,\,h\in H$ considered a set is a finite union of pairwise distinct right cosets $G h_i,\,h_i \in H$. Then the Hecke algebra $\mathcal{H} (G, H)$ over $\Z$ is defined to be the free $\Z$-module with basis $\{G h G : h \in H\}$. 
To define the multiplication of two double cosets $G h G$ and $G \tilde h G$ fix one system of representatives $h_i$ and $\tilde h_j$ for their right coset decompositions. There is a well defined multiset of double cosets $[G \hat h_i G]_i$ such that $[G h_i \tilde h_j]_{i, j} = \bigcup_i \bigcup_{j} [G \hat h_{i,j}]$ where $G \hat h_i G = \dot{\bigcup}_j G \hat h_{i,j}$. Hence the product $G h G \cdot G \tilde h G = \sum_i G \hat h_i G$ is well defined.

We will call the Hecke algebras associated to groups which have coefficients in $\Lambda$ the integral Hecke algebra. The terms local Hecke algebra and adelic Hecke algebra will be use analogously.

We will treat several settings in parallel. In. the unimodular case we assume $\Omega$ to be a finite central simple $K$ algebra.
In the modular case we distinguish the hermitian and the quaternionic case.
We either assume that $\Omega | K$ is an extension of fields of degree $2$ and $\bar \cdot$ is the Galois conjugation or that $\Omega$ is a quaternion algebra and $\bar \cdot$ is its canonical involution.

We first point out that the approximation lemmata given in \cite{R09} can be used to show that the right coset decomposition of any double coset and its adelic decomposition into right coset are in one to one correspondence.

\begin{lemma}Let $M \in \Inv{n} (\Lambda)$ and $\GL{n} (\adelic{\Lambda}) M \GL{n} (\adelic{\Lambda}) = \dot \bigcup_i \GL{n} (\adelic{\Lambda}) M_i$ with $M_i \in \Inv{n} (\adelic{\Lambda})$.
Then there are integral matrices $\widetilde{M}_i \in \Inv{n} (\Lambda)$ such that $\GL{n} (\adelic{\Lambda}) M_i = \GL{n} (\adelic{\Lambda}) \widetilde{M}_i$ and $\GL{n} (\Lambda) M \GL{n} (\Lambda) = \dot \bigcup_i \GL{n} (\Lambda) \widetilde{M}_i$.

An analogous statement holds true for $M \in \Delta_n (\Lambda, m),\, m \in \mathfrak{o}_K \setminus \{0\}$.
\end{lemma}
We only \textbf{Prove} the unimodular case, since the modular case is completely analogous.

Suppose there is an adelic right coset $\GL{n} (\adelic{\Lambda}) \hat M \subseteq \GL{n} (\adelic{\Lambda}) M \GL{n} (\adelic{\Lambda})$ with $\hat M \in \Inv{n} (\adelic{\Lambda})$.
Then there is a matrix $V \in \RGL{n} (\adelic{\Lambda})$ such that $\GL{n} (\adelic{\Lambda}) M V = \GL{n} (\adelic{\Lambda}) \hat M$.
Choose $m \in \mathfrak{o}_K \setminus \{0\}$ such that $m M^{-1},\, m {\hat M}^{-1} \in \Mat{n} (\adelic{\Lambda})$ and an approximation $\tilde V \in \GL{n} (\Lambda)$ modulo $m \adelic{\Lambda}$ of $V$.
Then $M \tilde V {\hat M}^{-1}$ and $\hat M (M \tilde V)^{-1}$ are integral in $\Mat{n} (\adelic{\Omega})$ and hence $\GL{n} (\adelic{\Lambda}) M \tilde V = \GL{n} (\adelic{\Lambda}) \hat M$.

On the other hand a standard argument shows that two right cosets which are distinct with respect to $\Lambda$ are also distinct with respect to $\adelic{\Lambda}$.\\
\vspace*{-5.5ex}\\
\hspace*{5cm}\hfill$ \Box$
\vspace{3ex}

This immediately yields a structure theorem for integral Hecke algebras.
\begin{theorem}
The integral unimodular and moduar Hecke algebras $\mathcal{H} (\GL{n} (\Lambda), \Inv{n} (\Lambda))$ and $\mathcal{H} (\Sp{n} (\Lambda), \Delta_n (\Lambda))$ are subalgebras of the adelic analogs $\mathcal{H} (\GL{n} (\adelic{\Lambda}), \Inv{n} (\adelic{\Lambda}))$ and $\mathcal{H} (\Sp{n} (\adelic{\Lambda}), \Delta_n (\adelic{\Lambda}))$, respectively. The embeddings map a double coset over $\Lambda$ generated by $M$ to the adelic double coset generated by $M$. 
\end{theorem}

We want to point out that the exceptional modular case $\Omega = K$ and $\bar \cdot = \mathrm{id}_K$, which corresponds to the symplectic group, builds up the complete picture. Here the integral and the adelic Hecke algebra are isomorphic.

The adelic Hecke algebra inherits a tensor product decomposition from $\adelic{\Lambda}$. So does the Hecke algebra over $\Lambda$ whenever it is a principle ideal domain. We will see that this is no longer true for the Hecke algebra over any $\Lambda$ which is not a PID. First evidence for this behavior in the hermitian modular case was found in \cite{Eb08} by means of elementary methods. This phenomena will be investigated in section \ref{sec:integralheckealgebra}.

\section{Adelic Hecke algebras}

The first question arising from the preceding section is how to understand the structure of the adelic Hecke algebras. This is easy for the unimodular case. We immediately see that
\begin{align*}
  \mathcal{H} (\GL{n} (\adelic{\Lambda}),\,\Inv{n} (\adelic{\Lambda}) )
\cong
  \bigotimes_{\mathfrak{p} \in \adelic{\Omega}^f}
  \mathcal{H} (\GL{n} (\compl{\Lambda}{p}),\, \Inv{n} (\compl{\Lambda}{p}) )
.
\end{align*}

The modular case is more involved. Surprisingly the quaternionic case is much easier since the canonical automorphism acts trivially on $\adelic{\Omega}^f$. This yields an analog decomposition of the adelic Hecke algebra as in the unimodular case.
\begin{align*}
  \mathcal{H} (\Sp{n} (\adelic{\Lambda}),\,\Delta_n(\adelic{\Lambda}))
\cong
  \bigotimes_{\mathfrak{p} \in \adelic{\Omega}^f}
  \mathcal{H} (\Sp{n} (\compl{\Lambda}{p}), \Inv{n} (\compl{\Lambda}{p}))
.
\end{align*}

One import case remains to be treated. This is the hermitian modular case. The Galois conjugation acts non trivially on split prime ideals in $\Lambda$. For the associated $\mathfrak{p} \in \adelic{\Omega}^f$ define the monoid of resticted invertible matrices $\RInv{n} (\compl{\Lambda}{p}) \subseteq \Inv{n} (\compl{\Lambda}{p}) \times \N_0$ to be the set of all pairs $(M, l)$ such that $\mathfrak{p}^l M^{-1}$ is integral.

Let us first present the result.
\begin{align*}
& \quad
  \mathcal{H} (\Sp{n} (\adelic{\Lambda}),\, \Delta_n (\adelic{\Lambda}) )
\\
&\cong
  \bigotimes_{\mathfrak {p} = \bar{\mathfrak{p}} \in \adelic{\Omega}^f}
  \mathcal{H} (\Sp{n} (\compl{\Lambda}{p}),\, \Delta_n (\compl{\Lambda}{p}) )
  \otimes
  \hspace{-2ex}
  \bigotimes_{\{\mathfrak{p},\, \bar{\mathfrak{p}}\} : \mathfrak{p} \ne \bar{\mathfrak{p}} \in \adelic{\Omega}^f}
  \hspace{-2ex}
  \mathcal{H} (\Sp{n} (\compl{\Lambda}{p} \oplus \ccompl{\Lambda}{p})),\, \Delta_n (\compl{\Lambda}{p} \oplus \ccompl{\Lambda}{p}) )
\\
&\cong
  \bigotimes_{\mathfrak {p} = \bar{\mathfrak{p}} \in \adelic{\Omega}^f}
  \mathcal{H} (\Sp{n} (\compl{\Lambda}{p}),\, \Delta_n (\compl{\Lambda}{p}) )
  \otimes
  \hspace{-2ex}
  \bigotimes_{\{\mathfrak{p},\, \bar{\mathfrak{p}}\} : \mathfrak{p} \ne \bar{\mathfrak{p}} \in \adelic{\Omega}^f}
  \hspace{-2ex}
  \mathcal{H} (\GL{n} (\compl{\Lambda}{p}),\,\RInv{n} (\compl{\Lambda}{p}) )
.
\end{align*}
Here $\GL{n} (\compl{\Lambda}{p})$ acts trivially on the second component of $\RInv{n} (\compl{\Lambda}{p})$. The first isomorphism is clear since the definition of modular matrices only involves $M$ and $\overline{M}$.

We have to \textbf{Prove} that 
\begin{align*}
  \mathcal{H} (\Sp{n} (\compl{\Lambda}{p} \oplus \ccompl{\Lambda}{p}),\, \Delta_n (\compl{\Lambda}{p} \oplus \ccompl{\Lambda}{p}) )
\cong
  \mathcal{H} (\GL{n} (\compl{\Lambda}{p}),\,\RInv{n} (\compl{\Lambda}{p}) )
.
\end{align*}
Write $(M_\mathfrak{p},\,M_{\bar{\mathfrak{p}}}) = M \in \Delta_n (\compl{\Lambda}{p} \oplus \ccompl{\Lambda}{p},\,m)$ where $m \in \mathrm{Fix}_{\bar \cdot} (\compl{\Lambda}{p} \oplus \ccompl{\Lambda}{p}) \setminus \{(0,0)\}$.
The restriction $M^* J M = m J$ is equivalent to $M_{\bar{\mathfrak{p}}} = m J^{-1} (M_\mathfrak{p} ^{-1})^* J$. Since $M_{\bar{\mathfrak{p}}}$ has to be integral $m M_\mathfrak{p} ^{-1}$ has to be integral, too, and this completely determines $M$. The right and double coset decomposition only depends on the first component and this proves the statement.\\
\vspace*{-5.5ex}\\
\hspace*{5cm}\hfill$ \Box$\\
\hspace*{3ex}

The theorems on generators of primary components given in \cite{Kr90} and \cite{Kr87_ha} generalize to all local Hecke algebras mentioned above.
To see this in the unimodular case remember $\compl{\Lambda}{p} = \Mat{r}(D)$ for some local divison algebra $D$ and some $r \in \N$. Then $\Inv{n} (\compl{\Lambda}{p}) \cong \Inv{r n} (D)$ and we can choose some prime element $\pi_D \in D$. The quaternionic modular case as well as the non split hermitian modular case behave completely analogously to the primary components of the Hecke algebra of the Hurwitz order and the PID hermitian Hecke algebras, respectively. In the split hermitian modular case we can proceed as in the unimodular case, since $\GL{2 n} (\compl{\Lambda}{p})$ leaves invariant the second component of $\RInv{2 n} (\compl{\Lambda}{p})$.

Let $[B_1,\ldots,B_i],\,i\in\N$ denote a block diagonal matrix.
\begin{lemma}
\begin{enumerate}[A)]
\item A set of generators of $\mathcal{H} (\GL{n} (\compl{\Lambda}{p}), \Inv{n} (\compl{\Lambda}{p}) )$ is given by the double cosets $T_i$ for $i=1,\ldots,r n$. 
Each $T_i$ is the double coset associated to $[I_{n r - i}, \pi_D I_i] \in \Mat{nr} (D) \cong \Mat{n} (\compl{\Lambda}{p})$.
Notice that for the Hurwitz order this corresponds to the generators satisfying $e_i \| e_{i+1}$ given in \cite{Kr87_ha}.
\item In the inert hermitian modular case suppose $\pi_D = \overline{\pi_D}$. A set of generators of $\mathcal{H} (\Sp{n} (\compl{\Lambda}{p}),\,\Delta_n (\compl{\Lambda}{p}))$ is given by $T_i$ with $i = 0,\ldots,r n$. 
Here $T_0$ is the double coset associated to $[I_{nr},\pi_D I_{nr}]$.
The $T_i,\,i=1,\ldots,r n$ are double cosets associated to $[I_{r n - i},\, \pi_D I_i,\pi_D^2 I_{r n - i},\,\pi_D I_i]$.
\item In the ramified quaternionic case and in the ramified hermitian modular case a set of generators of $\mathcal{H} (\Sp{n} (\compl{\Lambda}{p}),\,\Delta_n (\compl{\Lambda}{p}))$ is given by $T_i$ with $i = 0,\ldots,r n$. 
The $T_i,\,i=0,\ldots,r n$ are double cosets associated to $[I_{r n - i},\, \pi_D I_i,\pi_D \overline{\pi_D} I_{r n - i},\,\pi_D I_i]$.
\item In the split hermitian modular case a set of generators of the Hecke algebra $\mathcal{H} (\GL{2 n} (\compl{\Lambda}{p}),\, \RInv{2 n} (\compl{\Lambda}{p}))$ is given by $T_i,\,i=0,\ldots,2 n$. The $T_i$ are double cosets generated by $[I_{2n - i},\pi_D I_i]$. The second component of each $T_i$ is $1$.
\item In the split quaternionic modular case a set of generators of the Hecke algebra $\mathcal{H} (\Sp{n} (\compl{\Lambda}{p}),\,\Delta_n (\compl{\Lambda}{p}))$ is given  by $S,\,S'$ and $R_2,\ldots,R_{2 n}$. Here $S$ is generated by $[I_{2 n}, \pi_D I_{2 n}]$ and $S'$ is generated by $[I_{2 n - 1},\pi_D,\pi_D I_{2n - 2}, 1, \pi_D]$. The $R_{2 i}$ are generated by $[I_{2n - 2 i}, \pi_D I_{2 i},\pi_D^2 I_{2n - 2 i}, \pi_D I_{2 i}]$ and the $R_{2 i - 1}$ are generated by $[I_{2 n - 2 i + 1}, \pi_D I_{2 i - 1}, \pi_D ^2 I_{2 n - 2 i}, \pi_D, \pi_D ^2, \pi_D I_{2 i - 2}]$.
\end{enumerate}

All these Hecke algebras are commutative and all sets of generators are algebraically independent. 
\end{lemma}
\begin{eproof}
Proceed as described above.
The last statements follow from the existence of the canonical antiisomorphism $\cdot^*$ of the Hecke algebra and from counting the double cosets with fixed maximal elementary divisors or fixed similitude, respectively.
\end{eproof}
\vspace{2ex}

To investigate the structure and determinate generators of the integral Hecke algebra we need to prove one further lemma. It is a far generalization of well known lemmata one the multiplication of certain elements of the primary components, yet easy to prove.

For the rest of this section we will only be concerned with the local Hecke algebras.
We will identify double cosets $M$ with the ordered set of their defining elementary divisor valuations $\hat e (M) := (\hat e_1 (M) \le \cdots \le \hat e_{r n} (M))$ in the unimodular case or $(\hat e (M) := ((\hat e_1 (M) \le \cdots \le \hat e_{r n} (M)), l)$ where $m \sim \mathfrak{p}^l$ and $2 \mid l$, if $\mathfrak{p}$ is ramified, in the modular case.
The split hermitian modular case is an exception.
In this case we will identify the double cosets with $(\hat e (M) := (\hat e_1 (M) \le \cdots \le \hat e_{2 n} (M)), l)$ where $m \sim (\mathfrak{p} \bar{\mathfrak{p}})^l$.
Here $\hat e_i (M)$ is the valuation of the $i$-th elementary divisor $e_i (M)$.

Remember the $\mathfrak{p}$-rank of a double coset. It its the rank of any generator  of the double coset in $\compl{\Lambda}{p} / \mathfrak{p} \compl{\Lambda}{p}$.

An important ingredient to prove the preceding lemma was that for any double coset $M$ we have $T_i M = P + \sum_j R_j$, where $\hat e (P) = \hat e (T_i) + \hat e (M)$ and $R_j$ is a double coset which has $\mathfrak{p}$-rank less than the $\mathfrak{p}$-rank of $P$.
We are now going to generalize this for the product of arbitrary double cosets. To do so let us introduce an ordering on the double cosets.
\begin{align*}
  (\hat e(M), l) \le (\hat e(M'), l')
&:\Leftrightarrow
  (l < l')
  \vee
  (l = l' \wedge \hat e(M) < \hat e (M'))
\\
  \hat e(M) < \hat e (M')
&:\Leftrightarrow
  \exists i : \hat e_i (M) > \hat e_i (M')
             \wedge \forall j < i : \hat e_j (M) = \hat e_j (M')
.
\end{align*}

An induction on the double cosets with respect to this ordering yields
\begin{lemma}\label{dcosetmultiplication}
Given double cosets $M$ and $N$ we have $M N = P + \sum_j R_j$ where $\hat e(P) = \hat e(M) + \hat e(N)$ and $R_j$ are double cosets satisfying $R_j < P$.
\end{lemma}

\section{The integral Hecke algebra}\label{sec:integralheckealgebra}

In section \ref{sec:subheckealgebra} we have already shown that the integral Hecke algebras are subalgebras of the adelic ones. We will now give a set of generators and prove some results on it, that renders their computation possible. One special case which gave rise to the presented investigation will be given explicitely. Here we can even prove the minimality of a set of generators with respect to inclusion.

Let $\hat e (M)$ denote the tuple of all local elementary divisor valuations $(\hat e_\mathfrak{p} (M))_{\mathfrak{p} \in \adelic{\Omega}^f}$.
Let $G$ be the set of all double cosets $M$ in $\Inv{n} (\Lambda)$ or $\Delta_n (\Lambda)$ such that whenever $\hat e (M) = \hat e (N_1) + \hat e (N_2)$ for two double cosets $N_1,\,N_2$ over $\Lambda$ it follows that one of them is trivial. In \cite{Eb08} Ensenbach called some double cosets satisfying this property indecomposable and for the sake of uniform naming so will we. Nevertheless due to \ref{dcosetmultiplication} one should think of them as irreducible elements of the Hecke algebra.

Using the definition of $G$ and lemma \ref{dcosetmultiplication} we immediately see that $G$ generates the Hecke algebra as a $\Z$-algebra. This renders the indecomposable double coset important and we want to deduce an upper bound for the number of places with non trivial local elementary divisors as well as on the number of distinct local elementary divisors.

We may regard any local elementary divisor as an element of the Picard group. In \cite{R09} we proved that an adelic double coset is the image of an integral double coset if and only if the sum of its elementary divisors is zero in $\mathrm{Pic}(\Lambda)$. Denote the order of the Picard group by $h_\Lambda$ and its exponent by $e_\Lambda$. Furthermore consider the fundamental divisors of a matrix or double coset. Set $\hat f_i := \hat e_i - \hat e_{i-1}$ with $\hat e_0 := (0)_{\mathfrak{p} \in \adelic{\Omega}^f}$.

Suppose $h_\Lambda \ne 1$. Any double coset with $\hat f_{\mathfrak{p}, i} \ge e_\Lambda$ for some $\mathfrak{p}$ and $i$ is decomposable.
Consider a double coset $M$ and assume $\hat f_{\mathfrak{p}, i} < e_\Lambda$ for all indices. Suppose that there are at least $b_\Lambda := (h_\Lambda - 1)^{e_\Lambda - 1} + 1$ pairwise distinct and non conjugated places $\mathfrak{p} \in \adelic{\Omega}^f$ such that $\hat f_\mathfrak{p} (M)$ is not trivial. Then we may find a subset of $e_\Lambda$ places $\mathfrak{p}_j$ and indices $i_j$ such that in $\mathrm{Pic}(\Lambda)$ we have $\sum_j \mathfrak{p}^{\hat f_{\mathfrak{p}_j, i_j}} \cdot (n + 1 - i_j)  = 0$. Hence $M$ is decomposable.

This proves
\begin{lemma}
The set $G$ of indecomposable double cosets generate the Hecke algebra. Moreover there is an upper bound to the number of places where any irreducible double coset is non trivial and there is an upper bound for the valuation of the local elementary divisors.
\end{lemma}

Since the indecomposability of any double coset only depends on the valuation of its elementary divisors and the images of the local left ideals in the Picard group, there are in deed only finitely many cases we have to consider.
In the quaternionic modular case the Picard group's exponent is $2$ and hence the set of indecomposable double cosets is easy to derive.
In the hermitian modular case the set of indecomposable double cosets may becomes rather big. Using the following idea the computation is feasible.

First consider a decomposition into cyclic $p$-groups $\mathrm{Pic} \cong \bigoplus_{i=1} ^r (C_{p_i ^{\eta_i}}) ^{l_i}$ with pairwise distinct pairs $(p_i, \eta_i)$. Notice that any integral double coset vanishes in each component of the Picard group. A double coset is indecomposable if and only if it is indecomposible with respect to each of these components. Hence we may assume that $\mathrm{Pic} (\Lambda) \cong C_{p^\eta}$ for one prime $p$ and $\eta \in \N$.

Consider the multisets $[c_j]$ of elements in $\mathrm{Pic} (\Lambda)$ such that $\sum_i c_i = 0$ and for any submultiset $[\tilde c_i]$ we have $\sum_i c_i \ne 0$. Given one of these multisets consider the multisets of pairs of places and indices $(\mathfrak{p}_j,\,i_j)$ satisfying $\mathfrak{p}_j (n + 1 - i_j) = c_i$. These multisets correspond in an one-to-one manner to indecomposable double cosets $M$ by virtue of
\begin{align*}
\forall \mathfrak{p} \forall i \in \{1,\ldots,n\} :
  \#\{j : \mathfrak{p}_j = \mathfrak{p} \vee i_j = i\}
= 
  \hat f_{\mathfrak{p},\,i} (M)
.
\end{align*}

To illustrate this helpful devise let us write down one example. We will calculate the indecomposable double cosets in the hermitian modular case with class number $2$ for $n=2$.

So lets first consider those cosets which have prime similitude $p \in \mathfrak{o}_K$. If $p$ is inert the valuation of the elementary divisors of irreducible cosets are $(0,0,1,1),\, (0,1,1,2)$ and $(1,1,1,1)$. If $p$ ramifies there are indecomposable double cosets $(1,1,1,1)$ and $(0,0,2,2)$. Moreover if the ramified prime ideal above $p$ is principle, $(0,1,1,2)$ is indecomposable, otherwise $(0,2,2,4)$ is indecomposable.

If $p$ splits we have to distinguish two cases.
If it splits into principle ideals all indecomposable double cosets have similitude $p$ and the valuation of the local elementary divisors over $\compl{\Lambda}{p}$ are $(0,0,0,0),\,(0,0,0,1)\,(0,0,1,1),\,(0,1,1,1)$ and $(1,1,1,1)$.
If it splits into non principle ideals those indecomposable double cosets with similitude $p$ are $(0,0,0,0),\,(0,0,1,1)$ and $(1,1,1,1)$. In addition there are three which have similitude $p^2$.
They are $(0,1,1,2),\,(0,2,2,2)$ and $(0,0,0,2)$.

There remain those double cosets which have similitude $p_1 p_2$ with distinct primes $p_1$ and $p_2$.
They are the product of exactly two adelic double cosets $(0,1,1,2)$ if $p_i$ ramifies and $(0,1,1,1)$ or $(0,0,0,1)$ with similitude $p_i$ if $p_i$ splits.

Actually this is a minimal generating set with respect to inclusion of sets. This can easily be seen by the algebraic independence of the generators of the adelic Hecke algebra.

\subsection{The law of interchange for the Siegel $\Phi$-operator}
The application of these result to the theory of modular forms renders the following question important. One is interested in proving that the Siegel-$\Phi$ operator is surjective. Suppose $\Lambda$ is a principle ideal domain, $(\#\Lambda ^\times) \mid r \in \N$ and $n \ge 2$ if $\Lambda$ is commutative, $n \ge 3$ otherwise. To define
\begin{align*}
  \Phi_r :
  \Q \otimes \mathcal{H} (\Sp{n} (\Lambda),\,\Delta_n (\Lambda))
  \rightarrow
  \Q \otimes \mathcal{H} (\Sp{n-1} (\Lambda),\,\Delta_{n-1} (\Lambda))
\end{align*}
consider the underlying modules of linear combinations of right cosets. Fix a set of right coset representatives $h_i$ of a double coset $M$ and suppose
\begin{align*}
  h_i
=
  \left(\begin{matrix}
  A_i & 0^{(n-1,1)} & B_i & * \\
  *   & \alpha_i    & *    & * \\
  C_i & 0^{(n-1,1)} & D_i  & * \\
  0^{(1,n-1)} & 0   & 0^{(1,n-1)} & \delta_i
  \end{matrix}\right)
\end{align*}
Then $\Phi_r(M) := \sum_i \delta_i ^r \Sp{n} (\Lambda) \left(\begin{smallmatrix}A_i & B_i \\ C_i & D_i\end{smallmatrix}\right) \Sp{n} (\Lambda)$.

Since we identified the integral Hecke algebra with the adelic one, our final theorem can be proved by a standard counting argument.
\begin{theorem}
The map $\Phi_r$ is surjective except if $\Lambda$ is non commutative and $r = 2 n - 1$ or $r = 2n - 2$.
\end{theorem}
For a \textbf{Proof} consider Krieg's extensive studies of the Hurwitz order in \cite{Kr87_ha}. Note that in the quaternionic case he writes $\delta^{r / 2}$ instead of $\delta^r$.
All his arguments are still valid, but we have to replace $p$ by $\# (\compl{\Lambda}{p} / \mathfrak{p} \compl{\Lambda}{p})$. This doesn't change the underlying systems of linear relations for the images of the Hecke algebra's generators.\\
\vspace*{-5.5ex}\\
\hspace*{5cm} \hfill$ \Box$

\bibliography{bibliography}{}

\begin{thebibliography}{Rau09}

\bibitem[Ens08]{Eb08}
Marc Ensenbach.
\newblock {\em {Hecke-Algebren zu unimodularen und unit\"aren Matrixgruppen
  \"uber Dedekind-Ringen}}.
\newblock PhD thesis, RWTH Aachen University, 2008.

\bibitem[Kri87]{Kr87_ha}
Aloys Krieg.
\newblock {The Hecke-algebras related to the unimodular and modular group over
  the Hurwitz order of integral quaternions.}
\newblock {\em Proc. Indian Acad. Sci., Math. Sci. (Ramanujan Birth Centenary
  Volume)}, 97(1-3):201--229, 1987.

\bibitem[Kri90]{Kr90}
Aloys Krieg.
\newblock {Hecke algebras.}
\newblock {\em Mem. Am. Math. Soc.}, 435:158 p., 1990.

\bibitem[Rau09]{R09}
Martin Raum.
\newblock Elementary divisor theory for the modular group over hermitian fields
  and definite quaternion algebras.
\newblock {arXiv:0907.2762v1 [math.NT]}, 2009.

\end{thebibliography}
\bibliographystyle{alpha}

\end{document}